\documentstyle[11pt, amsfonts]{amsart}
\begin{document}
\def\R{{\mathbb R}}
\def\Rad{{\mathcal R}}
\def\Z{{\mathbb Z}}
\def\C{{\mathbb C}}
\newcommand{\trace}{\rm trace}
\newcommand{\Ex}{{\mathbb{E}}}
\newcommand{\Prob}{{\mathbb{P}}}
\newcommand{\E}{{\cal E}}
\newcommand{\F}{{\cal F}}
\newtheorem{df}{Definition}
\newtheorem{theorem}{Theorem}
\newtheorem{lemma}{Lemma}
\newtheorem{pr}{Proposition}
\newtheorem{co}{Corollary}
\newtheorem{pb}{Problem}
\def\n{\nu}
\def\sign{\mbox{ sign }}
\def\a{\alpha}
\def\N{{\mathbb N}}
\def\A{{\cal A}}
\def\L{{\cal L}}
\def\X{{\cal X}}
\def\F{{\cal F}}
\def\c{\bar{c}}
\def\v{\nu}
\def\d{\delta}
\def\diam{\mbox{\rm dim}}
\def\vol{\mbox{\rm Vol}}
\def\b{\beta}
\def\t{\theta}
\def\l{\lambda}
\def\e{\varepsilon}
\def\colon{{:}\;}
\def\pf{\noindent {\bf Proof :  \  }}
\def\endpf{ \begin{flushright}
$ \Box $ \\
\end{flushright}}

\title[Measure comparison and distance inequalities]
{Measure comparison and distance inequalities for convex bodies}

\author{Alexander  Koldobsky}
\address[Alexander  Koldobsky]{Department of Mathematics, University of Missouri, Columbia, MO 65211, USA}
\email{koldobskiya@@missouri.edu}
\author{Grigoris Paouris}
\address[Grigoris Paouris]{Department of Mathematics, Texas A\&M University, College Station, TX 77843, USA}
\email{grigoris@@math.tamu.edu}
\author{Artem Zvavitch}
\address[Artem Zvavitch]{Department of Mathematical Sciences\\ Kent State University\\ Kent, OH USA}
\email{zvavitch@@math.kent.edu}
\thanks{The first named author was supported by the U.S. National Science Foundation Grant DMS-1700036.
He also thanks the Max Planck Institute for Mathematics for support during his stay there in May-June 2019. 
The second named author was supported by the NSF grant DMS-1812240.
The third author was supported in part by  the U.S. National Science Foundation Grant DMS-1101636, Simons Foundation and by La Comue Universit\'e Paris-Est and B\'ezout Labex of Universit\'e Paris-Est.}
\date{}

\begin{abstract}  We prove new versions of the isomorphic Busemann-Petty problem for two different
measures and show how these results can be used to recover slicing and distance inequalities.  We also prove a sharp upper estimate for the outer volume ratio distance from an arbitrary convex body to the unit balls of subspaces of $L_p$.

\end{abstract}  
\maketitle
\section{Introduction}

Let $f$ be a non-negative locally integrable function on $\R^n$, and let $\mu$ be the measure in $\R^n$ with density
$f$, i.e. for every compact set $B\subset \R^n$
$$
\mu(B)=\int\limits_B f(x)dx.
$$
For $\xi\in S^{n-1},$ let $\xi^\bot$ be the central hyperplane
orthogonal to $\xi$, i.e. $\xi^\bot=\{x\in \R^n:\ \langle x,\xi\rangle=0\}.$ 
We write
$$\mu(B\cap \xi^\bot)=\int_{B\cap \xi^\bot} f(x)\ dx,$$
where in the right-hand side we integrate the restriction of $f$ to $\xi^\bot$ with respect to
 Lebesgue measure on $\xi^\bot.$

The following Busemann-Petty problem for general measures
(BPGM) was solved in \cite{Zv}. Fix $n\ge 2$. Given two
convex origin-symmetric bodies $K$ and $L$ in $\R^n$ such
that
$$
\mu(K\cap \xi^\bot)\le\mu(L\cap \xi^\bot)
$$
for every $\xi \in S^{n-1}$, does it follow that
$$
\mu(K)\le \mu(L)?
$$
It was proved in \cite{Zv}, that for every strictly positive even continuous density $f$ the answer to BPGM is
 affirmative if $n \le 4$ and negative if $n \ge 5$.
 
The BPGM problem is a generalization of the original Busemann-Petty  problem,
posed in 1956 in \cite{BP} and   asking the same question for 
Lebesgue measure $\mu(K)=|K|$, where the density $f\equiv 1$; see \cite{Z1, GKS, G3, K1}  
for the last steps of the solution and historical details. Throughout the paper we denote by $|K|$
the volume of $K$ in proper dimension.
 
Since the answer to BPGM is negative in most dimensions, it is natural to consider the following isomorphic version.
Let $C_n$ be the smallest constant  $C$ such that for any measure $\mu$ with continuous non-negative even density $f$ and any  two origin-symmetric convex bodies $K$ and $L$ in $\R^n$ satisfying
$$
\mu(K\cap \xi^\bot)\le\mu(L\cap \xi^\bot)
$$
for every $\xi \in S^{n-1}$, one necessarily has
\begin{equation}\label{izom}
\mu(K)\le  C \mu(L).
\end{equation}
It was shown in \cite{KZ} that $C_n\le  \sqrt{n}.$ It is not known 
whether this estimate is optimal.   The above result is a particular case of a more general statement connecting the BPGM with the Banach-Mazur distance between two convex bodies. Let
$$d_{BM}(K, {\cal{I}}_n)=\inf\{\frac{b}{a}:\ a D\subseteq K\subseteq b D, \mbox{ for some } D \in {\cal I}_n \}$$ be the Banach-Mazur distance between symmetric bodies $K \subset \R^n$ and the class of intersection bodies (see equality  (\ref{genint}), below, for the definition). It was shown in   \cite{KZ} that the constant in (\ref{izom}) is actually less or equal to $d_{BM}(K, {\cal{I}}_n)$. It is an open question whether the Banach-Mazur distance can be replaced here by the outer volume ratio distance $d_{\rm ovr}(K, {\cal{I}}_n);$ see the definition below. This would improve the $\sqrt{n}$ estimate for certain classes
of bodies, for example we could get an absolute constant for unconditional convex bodies; see \cite{K2}.

In this note we consider a version of the isomorphic Busemann-Petty problem with two different measures,
as follows. Define the outer volume ratio distance from a star body $K$ in $\R^n$ to the
class of generalized $k$-intersection bodies (see equality  (\ref{genint}), below, for the definition) by
$$d_{\rm ovr}(K,{\cal{BP}}_k^n)= \inf \left\{\left(\frac{|D|}{|K|}\right)^{\frac 1n}:\ K\subset D,\ D\in {\cal{BP}}_k^n\right\}.$$

We prove
\begin{theorem} \label{bp-sections-1} Let $K,L$ be star bodies in $\R^n,$ let $0<k<n,$
and let $\mu_1,\ \mu_2$ be measures on $\R^n$ with non-negative locally integrable densities $f, g$ so that $\|g\|_\infty=g(0)=1$ and
\begin{equation}\label{bp1}
\mu_1(K\cap H)\le \mu_2(L\cap H),\qquad \forall H\in Gr_{n-k}.
\end{equation}
Then
$$\mu_1(K)\le (d_{\rm ovr}(K,BP_k^n))^k \frac n{n-k} |K|^{\frac kn}\left(\mu_2(L)\right)^{\frac{n-k}n}.$$
\end{theorem}

Note that by John's theorem \cite{J}, if $K$ is origin-symmetric convex then $d_{\rm ovr}(K,BP_k^n)\le \sqrt{n}.$
Also the distance $d_{\rm ovr}(K,BP_k^n)$ is smaller than  the Banach-Mazur distance.  In particular,  $d_{\rm ovr}(K,BP_k^n)$ is bounded by an absolute constant for many classes of convex
bodies $K$. For example, for unconditional convex bodies $K$ this distance is bounded by $\sqrt{e}$; see \cite{K2}.

A version of the Busemann-Petty problem for the moments of a measure on convex bodies
was established in \cite{BKK}. Similarly to the case of sections, we extend this result to the case
of two different measures.

For $p> 0$ denote by $L_p^n$ the class of the unit balls of $n$-dimensional 
subpaces of $L_p.$ In other words (see for example \cite[p.117]{K1}), $L_p^n$ is the class of origin-symmetric
convex bodies $D$ in $\R^n$ such that there exists a finite Borel measure $\nu_D$ on $S^{n-1}$
satisfying
\begin{equation}\label{subspLp}
\|x\|_D^p=\int_{S^{n-1}} |\langle x,\theta\rangle |^p\ d\nu_D(\theta),\qquad \forall x\in \R^n.
\end{equation}
For a convex body $K$ in $\R^n$, denote by
\begin{equation}\label{p-ovr}d_{\rm {ovr}}(K,L_p^n) = \inf \left\{ \left( \frac {|D|}{|K|}\right)^{1/n}:
\ K\subset D,\ D\in L_p^n \right\}\end{equation}
the outer volume ratio distance from $K$ to the class of unit balls of subspaces of $L_p.$
\medbreak

\begin{theorem}\label{bp-g} Let $K,M \subset \R^n$ be two star bodies, let $p>0,$
and let $f,g$ be two non-negative measurable functions on $\R^n$, such that $\|g\|_\infty =g(0)=1.$
Suppose that for every $\xi\in S^{n-1}$
\begin{equation}\label{bpg}
\int_K|\langle x,\xi\rangle |^p g(x) dx \le \int_M|\langle x,\xi\rangle |^p f(x)\ dx.
\end{equation}
Then
$$
  \left(\int_{K}{ g}(x) dx \right)^{(n+p)/n} \le  \frac {n+p}n  d_{\rm ovr}(M,L_p^n)^p\ |M|^{\frac{p}{n}} \int_M  f(x)dx.
$$
\end{theorem}

We show how Theorems \ref{bp-sections-1} and \ref{bp-g} can be used to recover slicing inequalities
for sections and moments from \cite{K2} and \cite{BKK} and distance inequalities from \cite{KK}
and \cite{BKK,KLi}. In particular, it was shown in \cite{BKK,KLi} that
$$c\sqrt{\frac np}\le \sup_K\ d_{\rm ovr}(K,L_p^n)\le \sqrt{n}.$$
We prove a sharp estimate
$$\sup_K\ d_{\rm ovr}(K,L_p^n)\le C\sqrt{\frac{n+p}p}.$$
Here $c,C$ are absolute constants.

\section{Sections}\label{sec:sec}

We need several definitions and facts.
A closed bounded set $K$ in $\R^n$ is called a {\it star body}  if 
every straight line passing through the origin crosses the boundary of $K$ 
at exactly two points different from the origin, the origin is an interior point of $K,$
and the {Minkowski functional} 
of $K$ defined by 
$$\|x\|_K = \min\{a\ge 0:\ x\in aK\}$$
is a continuous function on $\R^n.$  We will denote by   ${\cal{K}}_{n}$ the class of convex bodies in $\R^n$. 
By $|K|_m$, or simply $|K|$ when there is no ambiguity, we denote the $m$-dimensional Lebesgue measure of a set $K$.
We use the polar formula for the volume $|K|$ of a star body $K:$
\begin{equation}\label{polar-vol}
|K|=\frac 1n \int_{S^{n-1}} \|\theta\|_K^{-n} d\theta.
\end{equation}

If $\mu$ is a measure on $\R^n$ with continuous density $f,$ we have
\begin{equation}\label{polar-meas}
\mu(K)=\int_K f = \int_{S^{n-1}} \left(\int_0^{\|\theta\|_K^{-1}} r^{n-1}f(r\theta) dr \right) d\theta.
\end{equation}

Let $Gr_{n-k}$ be the Grassmanian of $(n-k)$-dimensional subspaces of $\R^n.$
For $1\le k \le n-1,$  the {$(n-k)$-dimensional spherical Radon transform} 
$\Rad_{n-k}:C(S^{n-1})\to C(Gr_{n-k})$ is a linear operator defined by
$$\Rad_{n-k}g (H)=\int_{S^{n-1}\cap H} g(x)\ dx,\quad \forall  H\in Gr_{n-k}$$
for every function $g\in C(S^{n-1}).$ 
\smallbreak
For every $H\in Gr_{n-k},$ the $(n-k)$-dimensional volume of the section of a star body $K$ by $H$
can be written as
\begin{equation}\label{vol-sect}
|K\cap H| = \frac1{n-k} \Rad_{n-k}(\|\cdot\|_K^{-n+k})(H).
\end{equation}
More generally, for a measure $\mu$ with continuous density $f$ and any $H\in Gr_{n-k}$, we write
\begin{equation}\label{meas-sect}
\mu(K\cap H)=\int_{K\cap H} f = \Rad_{n-k}\left(\int_0^{\|\cdot\|_K^{-1}} r^{n-k-1}f(r\ \cdot)\ dr \right)(H).
\end{equation}

A generalization of the concept of an intersection body was introduced by Zhang \cite{Z2}. 
We say that an origin symmetric star body $D$ in $\R^n$ is a {generalized $k$-intersection body}, 
and write $D\in {\cal{BP}}_k^n,$  if there exists a finite Borel non-negative measure $\nu_D$
on $Gr_{n-k}$ so that for every $g\in C(S^{n-1})$
\begin{equation}\label{genint}
\int_{S^{n-1}} \|x\|_D^{-k} g(x)\ dx=\int_{Gr_{n-k}} R_{n-k}g(H)\ d\nu_D(H).
\end{equation}
When $k=1$ we get the original class of intersection bodies ${\cal{BP}}_1^n={\cal{I}}_n$ introduced by Lutwak \cite{L2}.
\medbreak

\noindent {\bf Proof of Theorem \ref{bp-sections-1}.} For a small $\delta>0,$ let $D\in {\cal{BP}}_k^n$ be a body such that $K\subset D$ and
\begin{equation}\label{sect11}
|D|^{\frac 1n}\le (1+\delta)\ d_{\rm ovr}(K,BP_k^n)\ |K|^{\frac 1n},
\end{equation} 
and let $\nu_D$ be the measure on $Gr_{n-k}$ corresponding to $D$ by the definition (\ref{genint}).

By (\ref{vol-sect}) and (\ref{meas-sect}), the condition (\ref{bp1}) of the theorem can be written as 
$$\Rad_{n-k}\left(\int_0^{\|\cdot\|_K^{-1}} r^{n-k-1}f(r\ \cdot)\ dr \right)(H) \le  
\Rad_{n-k}\left(\int_0^{\|\cdot\|_L^{-1}} r^{n-k-1}g(r\ \cdot)\ dr \right)(H)$$
for every $H\in Gr_{n-k}.$
Integrating both sides of the latter inequality with respect to $\nu_D$ and using the definition (\ref{genint}), we get
\begin{equation}\label{integration}
\int_{S^{n-1}} \|x\|_D^{-k} \left(\int_0^{\|x\|_K^{-1}} r^{n-k-1}f(rx)\ dr \right)dx 
\end{equation}
$$\le  \int_{S^{n-1}} \|x\|_D^{-k} \left(\int_0^{\|x\|_L^{-1}} r^{n-k-1}g(rx)\ dr \right)dx,$$
which is equivalent to
\begin{equation}\label{sect12}
\int_K \|x\|_D^{-k}f(x)dx \le \int_L \|x\|_D^{-k}g(x)dx.
\end{equation}

Since $K\subset D,$ we have $1\ge \|x\|_K\ge \|x\|_D$ for every $x\in K.$ Therefore,  
$$\int_K \|x\|_D^{-k}f(x) dx \ge \int_K \|x\|_K^{-k}f(x) dx \ge \int_K f(x) dx=\mu_1(K).$$

On the other hand, by the Lemma from section 2.1 from Milman-Pajor
\cite[p.76]{MP}, 
$$
\left(\frac{\int_{L}\|x\|_D^{-k} {g}(x) dx}{ \int_{D}\|x\|_D^{-k} dx} \right)^{1/(n-k)}  \le \left(\frac{\int_{L}{g}(x) dx \large}{  \int_{D} dx} \right)^{1/n}.
$$
Since $\int_D\|x\|_D^{-k} dx =\frac{n}{n-k} |D|,$  (see \cite{MP}) we can estimate the right-hand side of
(\ref{sect12}) by
$$\int_L \|x\|_D^{-k} g(x) dx \le \frac n{n-k}  (\mu_2(L))^{\frac {n-k}n} |D|^{\frac kn}.$$
Now apply (\ref{sect11}) and send $\delta$ to zero to get the result. \endpf

In the case where $\mu_2$ is volume, i.e. $g\equiv 1,$ we get the following.

\begin{co} \label{bp-sections} Let $K,L$ be star bodies in $\R^n,$ let $0<k<n,$
and let $\mu$ be a measure on $\R^n$ with density $f$ so that
\begin{equation}\label{bp}
\mu(K\cap H)\le |L\cap H|,\qquad \forall H\in Gr_{n-k}.
\end{equation}
Then
$$\mu(K)\le (d_{\rm ovr}(K,BP_k^n))^k \frac n{n-k} |L|^{\frac {n-k}n} |K|^{\frac kn}.$$
\end{co}

Corollary \ref{bp-sections} provides a new proof of the slicing inequality for measures established earlier in \cite{K2}, as follows.
Let $c_{n,k}= |B_2^n|^{\frac {n-k}n}/|B_2^{n-k}|,$ where $B_2^n$ is the unit Euclidean ball in $\R^n.$
Note that $c_{n,k}\in (e^{-k/2},1)$ (see for example \cite[Lemma 2.1]{KL}).

\begin{co} \label{lowdim} (\cite{K2}) Let $K$ be a star body in $\R^n.$ Then for any measure $\mu$
with density $f$ on $K$ we have
$$\mu(K)\le \left(d_{\rm ovr}(K,{{\cal{BP}}_k^n})\right)^k\ \frac n{n-k} 
c_{n,k} \max_{H\in Gr_{n-k}} \mu(K\cap H)\ |K|^{k/n}.$$
\end{co}

\pf Apply Corollary \ref{bp-sections} to the bodies $K$ and $L=cB_2^n,$ where 
$$c=\left(\frac{\max_{H\in Gr_{n-k}} \mu(K\cap H)}{|B_2^{n-k}|}\right)^{\frac 1{n-k}}.$$
\endpf

We note that using a similar argument one can provide a sharp estimate for the constant in the  following isomorphic version of an  Busemann-Petty problem

\begin{co} \label{bp-sections-isom} Consider a constant $C_n>0$ such that for any convex, symmetric bodies $K,L$  in $\R^n,$ 
and any  $\mu$ be a measure on $\R^n$ with even density $f$ if
$$
\mu(K\cap \xi^\perp)\le |L\cap \xi^\perp|,\qquad \forall  \xi \in S^{n-1},
$$
one necessarily has
$$\mu(K)\le C_n |L|^{\frac {n-1}n} |K|^{\frac 1n}.$$  Then $c_1 \sqrt{n} \le C_n \le c_2 \sqrt{n}$, where $c_1, c_2$ are positive absolute constants.
\end{co}

\pf  The fact that $C_n \le c_2 \sqrt{n}$ follows immediately from $d_{\rm ovr}(K, BP_1^n) \le \sqrt{n}$ and  Corollary \ref{bp-sections}. To prove the lower bound on $C_n$ one can use an argument similar to that in the proof of Corollary \ref{lowdim} and use the bounds provided in \cite{KK, KLi}.
\endpf

As it was mentioned in the proof of Corollary  \ref{bp-sections-isom}, in the case $k=1,$ where ${\cal{BP}}_1^n={\cal{I}}_n$ is the class of intersection bodies, it was proved
in \cite{KK} (with an extra logarithmic term, which can be removed due to the result from \cite{KLi}) that
there exists an absolute constant $c>0$ so that
\begin{equation}\label{In-1}
c\sqrt{n}\le \max_{K \in  {\cal{K}}_{n}} d_{\rm ovr}(K,{{\cal{I}}_n})\le \sqrt{n},
\end{equation}
where maximum is taken over all origin-symmetric convex bodies, and the upper 
bound follows from John's theorem \cite{J}, see also \cite{AGM}.

On the other hand, for many classes of bodies, the outer volume ratio distance $d_{\rm ovr}(K,{{\cal{BP}}_k^n})$ 
admits better bounds. For example, this distance is bounded by an absolute constant for unconditional convex bodies \cite{K2}, and for the
unit balls  of subspaces of $L_p$ it is bounded by $\sqrt{p}$; see \cite{M1, KP}. 

A weaker version of the left-hand side inequality in \eqref{In-1} has been proved in [\cite{KPZ}, Proposition 6.1] where it was shown that 
$$ \max_{K \in  {\cal{K}}_{n}} d_{BM}(K,{{\cal{I}}_n}) \geq  d_{BM}( B_{\infty}^{n} , {\cal{I}}_{n}) \geq c\sqrt{n}, $$ 
where $d_{BM}$ is the Banach-Mazur distance. Actually the same proof provides bounds for the volume ratio, that complements \eqref{In-1}. In particular we have the following
\begin{equation}
\label{In-2}
c\sqrt{n}\le \max_{K \in  {\cal{K}}_{n}} d_{\rm vr}(K,{{\cal{I}}_n})\le \sqrt{n}.
\end{equation}
Here
$$d_{\rm vr}(K,{\cal{I}}_n)= \inf \left\{\left(\frac{|K|}{|D|}\right)^{\frac 1n}:\ D\subset K,\ D\in {\cal{I}}_n\right\}.$$

For completeness we briefly sketch the proof.
Recall that the volume ratio ${\rm v.r.}(K, L)$ of two convex bodies $K$ and $L$ is defined by 
${\rm v.r.}(K, L) = \inf\{(\frac{|K|}{|T(L)|})^{1/n}\}$, where the infimum is over all affine transformations $T$ of 
$\R^n$ for which $T(L)\subset K.$ We write ${\rm v.r.}(K,B_2^n)={\rm v.r.}(K).$ 
Let $ {\cal{VR}}_{n}(a)$ be the class of all symmetric convex bodies in $\R^{n}$ that have volume ratio less than $a$. Then, there exists a universal constant $c$ such that for every $n$, $ {\cal{I}}_{n} \subseteq  {\cal{VR}}_{n}(c)$ (see Proposition 6.2 in \cite{KPZ}). So, $d_{\rm vr}(K,{{\cal{I}}_n}) \geq d_{\rm vr}(K,{{\cal{VR}}_n}(c))$. By the sub-additivity of $ \log({\rm v.r.}),$ we get that for all $K_{1}, K_{2}, K_{3}$ centrally symmetric convex bodies in $\R^n$, 
$$ {\rm v.r.}( K_{1}, K_{2}) \leq {\rm v.r.}( K_{1}, K_{3})\ {\rm v.r.}( K_{3}, K_{2}) .$$
We choose $K_{1}:= B_{\infty}^{n}$, $K_{2}:= B_{2}^{n}$ and $ K_{3} \in {{\cal{VR}}_n}(c)$ and we get that 
$ d_{\rm vr}(B_{\infty}^{n},{{\cal{VR}}_n}(c)) \geq \frac{c^{\prime} \sqrt{n}}{ c}$. This proves the left-hand side inequality in \eqref{In-2}. The right-hand side follows from John's theorem \cite{J}.

It was proved in \cite{KPZ} that for any origin-symmetric convex body $K$ in $\R^n$ and for all $1\le k<n,$
$$d_{\rm ovr}(K,{{\cal{BP}}_k^n})\le \sqrt{\frac nk} \left(\log\left(\frac{en}k\right)\right)^{3/2}.$$

\begin{pb} Does there exist an absolute constant $c>0$ so that for all $2\le k<n$
$$c\sqrt{\frac nk}\le \max_{K \in  {\cal{K}}_{n}} d_{\rm ovr}(K,{{\cal{BP}}_k^n}),$$
where maximum is taken over all origin-symmetric convex bodies $K?$
\end{pb}
Actually, a weaker version of the above question was asked in \cite{KPZ} with $d_{\rm ovr}$ replaced by 
$d_{BM}.$ This question is also open.

\section{Moments}

\noindent {\bf Proof of Theorem \ref{bp-g}.} For $\delta>0,$ let $D$ be an origin-symmetric convex body 
such that $D\in L_p^n,\ M\subset D,$
and
\begin{equation} \label{ovr-cond}
|D|^{\frac 1n}\le (1+\delta)\ d_{\rm ovr}(M,L_p^n)\ |M|^{\frac 1n}.
\end{equation}
Let $\nu_D$ be the measure on $S^{n-1}$ corresponding to $D$ by (\ref{ovr-cond}).
Integrating both sides of (\ref{bp}) over $S^{n-1}$ with the measure $\nu_D$ we get
\begin{equation}\label{D-ineq}
\int_K\|x\|_D^p\ g(x) \ dx \le \int_M\|x\|_D^p\ f(x)\ dx.
\end{equation}
Since $M\subset D,$ we have $\|x\|_D\le \|x\|_M\le 1$ for every $x\in M,$ so the right-hand side of (\ref{D-ineq})
can be estimated from above by
$$\int_M\|x\|_D^p\ f(x) dx\le \int_M f(x)\ dx.$$
To estimate the left-hand side of (\ref{D-ineq}) from below, we use the Lemma from section 2.1 from Milman-Pajor
\cite[p.76]{MP}, which asserts  that if ${\tilde g}:\R^n \to \R$ a measurable function such that $\|{\tilde g}\|_\infty=1$ and $D$ is a symmetric convex body in $\R^n$, then 
$$
\left(\frac{\int_{\R^n}\|x\|_D^p {\tilde g}(x) dx}{ \int_{D}\|x\|_D^p dx} \right)^{1/(n+p)}  \ge \left(\frac{\int_{\R^n}{\tilde g}(x) dx \large}{  \int_{D} dx} \right)^{1/n}.
$$
Using that $\int_D\|x\|_D^p dx =\frac{n}{n+p} |D|$ we get
$$
\int_{\R^n}\|x\|_D^p {\tilde g}(x) dx \ge \frac{n}{n+p}|D|^{-p/n} \left(\int_{\R^n}{\tilde g}(x) dx \right)^{(n+p)/n}.
$$

applying this inequality for ${\tilde g}= g\cdot \chi_K$ we get

$$
\int_K\|x\|_D^p\ g(x) dx \ge \frac n{n+p} |D|^{-p/n}  \left(\int_{K}{ g}(x) dx \right)^{(n+p)/n}.
$$
Combining these inequalities with (\ref{ovr-cond}), we get
$$
  \left(\int_{K}{ g}(x) dx \right)^{(n+p)/n} \le   \frac {n+p}n |D|^{p/n}  \int_M   \|x\|_D^p f(x)\ dx \
 $$
$$
  \left(\int_{K}{ g}(x) dx \right)^{(n+p)/n} \le  \frac {n+p}n  d_{\rm ovr}(M,L_p^n)^p\ |M|^{\frac{p}{n}} \int_M   \|x\|_D^p f(x)dx
$$
\endpf

Taking $g\equiv 1$ we get the following corollary of Theorem \ref{bp-g}:

\begin{co}\label{bp-volume} Let $K,M$ be origin-symmetric convex bodies in $\R^n,$ let $p>0,$
and let $f$ be a non-negative continuous function on $\R^n.$
Suppose that for every $\xi\in S^{n-1}$
\begin{equation}\label{bp}
\int_K|\langle x,\xi\rangle|^p dx \le \int_M|\langle x,\xi\rangle |^p f(x)\ dx.
\end{equation}
Then
$$|K|^{\frac {n+p}n}\le \frac{n+p}n \left(d_{\rm {ovr}}(M,L_p^n)\right)^{p} |M|^{\frac pn}\int_M f.$$
\end{co}

\section{The distance to subspaces of  $L_p$}

Similarly to methods used in Section \ref{sec:sec}, choosing $K$ in Corollary  \ref{bp-volume} as a  proper multiple of the Euclidean ball, one gets the following inequality established earlier in \cite{BKK}: there exists an absolute constant $C$ so that for every origin-symmetric convex body $M$ and every $p\ge 1,$
$$
\min_{\xi \in S^{n-1}} \int_M|\langle x,\xi\rangle |^p f(x)\ dx  \, \le \,
(Cp)^{p/2}\, d^p_{\rm {ovr}}(M,L_p^n)\ |K|^{p/n} \int_K f(x)\,dx.
$$
Using this inequality, it was shown in \cite{BKK, KLi} (in \cite{BKK} there was an extra logarithmic 
term in the left-hand side, which was removed in \cite{KLi}) that there exists an absolute constant $c>0$
so that for every $n$ and every $p\ge 1$
\begin{equation}\label{dovrLp}
c\sqrt{\frac np}\le \max_M\  d_{\rm {ovr}}(M,L_p^n) \le \sqrt{n},
\end{equation}
where maximum is taken over all origin-symmetric convex bodies $M$ in $\R^n.$

The  goal of this section is to provide an  optimal estimate in the right-hand side of (\ref{dovrLp}) with respect to $p$.
Actually we will show that the upper bound can be achieved by considering a mush smaller class of convex bodies. Recall that  $$
K^\circ =\{x \in \R^n:  \langle x, y \rangle \le 1, \mbox{  for all } y \in K \}
$$
is the polar body of $K$ (we refer to \cite{S2, AGM}  for more information). Let $ p\geq 1$ and $C$ be a symmetric convex body in $\R^{n}$ of volume $1$. We define a convex body  $ Z_{p}(C)$ via its  support function
\begin{equation}
\label{Zp-def}
h_{Z_{p}(C) } ( \theta) := \left( \int_{C} | \langle x, \theta \rangle |^{p} dx\right)^{\frac{1}{p}} , \theta \in S^{n-1},
\end{equation}
These bodies has been introduced by the second-named author in \cite{Pa-01}, \cite{Pa-02} in order to investigate concentration properties of the random vector $X$ that is uniformly distributed in $C$. For results of this type see \cite{Pa-1}, \cite{Pa-sb}, \cite{FGP}. The $Z_{p}(C)$-bodies is a family of convex bodies that interpolate (in an increasing matter) between a multiple of a ball $(p=2)$ and the body  $C$, $(p=\infty)$, with a no-more than linear ``speed" in the distance (See \cite{Pa-02}). A more detailed collection of properties of these bodies can be found in \cite{Pa-1}, \cite{Pa-sb} or in  
  Chapter 5 in \cite{BGVV}. We define $ {\cal{M}}_{p}^{n}$ be the set of symmetric convex bodies in $ \R^{n}$ that are polar to some $Z_{p}(C)$, i.e. 
$$  {\cal{M}}_{p}^{n} := \left\{ L\in {\cal{K}}_{n} : \exists C\in K_{n} , |C|=1 , L= Z_{p}^{\circ}(C) \right\} .$$
We have that 
$$  {\cal{M}}_{p}^{n} \subseteq L_{p}^{n} . $$
Indeed,   writing \eqref{Zp-def} in polar coordinates 
\begin{eqnarray*}
  \| \theta \|_{L}^p = h^p_{Z_{p}(C) } ( \theta) &=& \int_{S^{n-1}} | \langle x, \theta \rangle |^{p}  \frac{\delta_{n}}{ \|x\|_{C}^{n+p} }d\sigma(x)\\ &=&   \int_{S^{n-1}} | \langle x, \theta \rangle |^{p}  d\mu(x)  ,
\end{eqnarray*}
where $ d\mu(x) :=  \frac{\delta_{n} }{ \|x\|_{C}^{n+p} }d\sigma(x)$ is a finite Borel measure on $S^{n-1}$ and $ \delta_{n}$ a constant depending only on $n$. It follows that $(\R^{n},  \| \cdot \|_{L})$ can be isometrically embedded into $ L_{p}$ (see Lemma 6.4 in \cite {K1}). 

At this point we have to emphasize that ${\cal{M}}_{p}^{n} $ is a much smaller class than $ L_{p}^{n} $. For example we have that if $ L\in  {\cal{M}}_{p}^{n}$, $ L= Z_{p}^{\circ}( C) $, 
$$ d_{BM} ( L, B_{2}^{n} ) = d_{BM} ( Z_{p}(C) , B_{2}^{n} ) \leq c p , $$
 while 
 $$ \max_{K\in L_{p}^{n}} d_{BM} ( K, B_{2}^{n} )\geq d_{BM} ( B_{p}^{n} , B_{2}^{n} ) \geq n^{ \frac{1}{2} - \frac{1}{p}} . $$
 In the above we have used the following well known application of Brunn-Minkowski inequality, 
 \begin{equation}
 \label{incl}
  Z_{p} ( C) \subseteq c p Z_{2} ( C)  ,
  \end{equation}
  where $ c>0$ is an absolute constant (see, for example, \cite{BGVV}, Proposition 5.1.2). 

 Another way to measure how different the two classes are is to consider the volume ratio. We will need  the following result of Lutwak and Zhang \cite{LZ} stating that for every star  body $C$ of volume $1$ and every $ p>0$, 
\begin{equation}
\label{LZ} 
| Z_{p}^{\circ} (C) |^{\frac{1}{ n}} \leq |Z_{p}^{\circ} ( D_{n})|^{\frac{1}{n}} \leq c_{1} \frac{1}{n}\sqrt{\frac{n+p}{p}},   
\end{equation} 
where $ D_{n}=B_2^n/|B_2^n| $ is the Euclidean ball of volume one and $ c_{1}>0$ an absolute constant. So, using \eqref{incl} and \eqref{LZ} we get that if  $L= Z_{p}^{\circ}( C) \in {\cal{M}}_{p}^{n}$,
$$ {\rm v.r.} (L) \leq \left( \frac{ | Z_{p}^{\circ}( C) |}{ | \frac{1}{cp} Z_{2}^{\circ} ( C) | }\right)^{\frac{1}{n}} \leq c^{\prime} \min\{ \sqrt{p}, \sqrt{n} \}  L_{C}, $$
where $ L_{C} := \left( \frac{|Z_{2} ( C)|}{ | B_{2}^{n} | } \right)^{\frac{1}{n}} $ is the isotropic constant of $C$ (we refer to \cite{BGVV} for a definition and properties).

\noindent  Moreover, by a result of K. Ball \cite{Ba1} (see also \cite{LYZ1}, Theorem 3) we have that 
$$ \max_{K \in L_{p}^{n} } {\rm v.r.} ( K) = {\rm v.r.} ( B_{p}^{n} ) \simeq n^{\frac{1}{2}-\frac{1}{p}} . $$
We should also notice that for  $ p\simeq n $ the two classes are isomorphically the same and  (isomorphically) equal to $ {\cal{K}}_{n}$ (the class of centrally symmetric convex bodies in $ \R^{n}$). Indeed, for $ p\geq n$ an application of Brunn-Minkowski inequality (see Lemma 3.6 in \cite{Pa-sb}) implies that for every $ C\in {\cal{K}}_{n}$, 
\begin{equation}
\label{BM-1}
c_{0} C \subseteq Z_{p} ( C) \subseteq C , 
\end{equation}
where $ c_{0} >0$ is an absolute constant.

\begin{theorem}\label{prop:overp}
There exists $c>0$ such that for every $ n\geq 2$, $p\geq 1$, 
\begin{equation}
\label{result-1}
 \sup_{K\in {\cal{K}}_{n}} d_{\rm {ovr}} \left( K, L_{p}^n\right) \leq \sup_{K\in {\cal{K}}_{n}} d_{\rm {ovr}} \left( K, {\cal{M}}_{p}^n\right)  \leq c\sqrt{\frac{n+p}{p}}
\end{equation} 
\end{theorem}

\pf
Let $K\in {\cal{K}}_{n}$, and let  $ C:= K^{\circ}/ | K^{\circ}|^{\frac{1}{n}}$. Note that 
\begin{equation}
\label{eq-1}
 |C|=1 \ {\rm and } \ C^{\circ}= |K^{\circ}|^{\frac{1}{n}} K.
 \end{equation}
Let $L\in {\cal{M}}_{p}^{n}$ be defined as 
\begin{equation}
\label{eq-2}
 L:= |K^{\circ}|^{-\frac{1}{n}}Z_{p}^{\circ}(C) \ {\rm and } \ L^{\circ} = | K^{\circ} |^{\frac{1}{n}} Z_{p} (C) . 
 \end{equation}
Note that, since $ |C|=1$, \eqref{Zp-def} implies that $ Z_{p} (C) \subseteq C$. So, by \eqref{eq-1} and \eqref{eq-2}, 
$$ C^{\circ} \subseteq Z_{p}^{\circ} ( C) \  {\rm or }  \  \frac{C^{\circ}}{ | K^{\circ}|^{\frac{1}{n}} } \subseteq \frac{Z_{p}^{\circ} ( C)}{ | K^{\circ}|^{\frac{1}{n}} } \  {\rm or } \ K \subseteq L. $$
We will also need the Bourgain-Milman theorem (\cite{BM}, see also \cite{AGM}) on lower bounds of the volume product of convex bodies. It states that that for every $C\subseteq  \R^{n}$ symmetric convex body of volume $1$, 
\begin{equation}
\label{BM}
| C^{\circ} |^{\frac{1}{ n}} \geq \frac{c_{0}}{ n}, 
\end{equation}
where $ c_{0}>0$ is an absolute constant. 
The last ingredient of the proof is an inequality of Lutwak and Zhang \cite{LZ} stating that for every star  body $C$ of volume $1$ and every $ p>0$, 
\begin{equation}
\label{LZ} 
| Z_{p}^{\circ} (C) |^{\frac{1}{ n}} \leq |Z_{p}^{\circ} ( D_{n})|^{\frac{1}{n}} \leq c_{1} \frac{1}{n}\sqrt{\frac{n+p}{p}},   
\end{equation} 
where $ D_{n}=B_2^n/|B_2^n| $ is the Euclidean ball of volume one and $ c_{1}>0$ an absolute constant. Using \eqref{eq-1}, \eqref{eq-2}, \eqref{BM} and \eqref{LZ}, we conclude that 
$$ \left( \frac{ | L | }{ | K |} \right)^{\frac{1}{ n}} = \left( \frac{ | Z_{p}^{\circ}  (C) |}{ | C^{\circ} |} \right)^{\frac{1}{n}} \leq c\sqrt{\frac{n+p}{p}} ,$$
where $c:= \frac{c_{1}}{ c_{0}}>0$ is an absolute constant. 
\endpf

\noindent  Note that using \eqref{BM-1}, \eqref{result-1} and \eqref{dovrLp} we have that for all $ p\geq 1 $, $ n\geq 2$, 
$$  \sup_{K\in {\cal{K}}_{n}} d_{\rm {ovr}} \left( K, L_{p}^n\right) \simeq \sqrt{\frac{n+p}{p}} . $$


\smallskip 

\noindent In the case where $ 1\leq p \leq \sqrt{n}$ there is an improvement in the inequality of Lutwak and Zhang due to Klartag and Milman \cite{KM}. Their result states that for every $ C$ symmetric convex body in $ \R^{n}$ of volume $1$ and every $ 1\leq p \leq \sqrt{n}$, 
\begin{equation}
\label{KM}
| Z_{p}^{\circ} (C) |^{\frac{1}{ n}} \leq   \frac{c_{2}}{nL_{C} }\sqrt{\frac{n}{p}},   
\end{equation}
where $c_{2}>0$ is an absolute constant and $L_{C}$ is the isotropic constant of $C$. If in the previous proof we will use \eqref{KM} instead of \eqref{LZ} we get the following

\begin{theorem}
\label{cor}
\noindent Let $n\geq 1$, $ K\in {\cal{K}}_{n}$ and $ 1\leq p \leq \sqrt{n}$. Then 
\begin{equation}
\label{result-2}
 d_{\rm {ovr}} \left( K, L_{p}^n\right) \leq \frac{c_{3}}{ L_{K^{\circ}}} \sqrt{\frac{n}{p}}
 \end{equation}
\end{theorem}
Let us also mention that one may define $M_{p}^{n}$ as the class of the polar bodies of $Z_{p}(\mu)$-bodies, where $ \mu$ is a log-concave measure. This class is isomorphically equivalent to the class that we used, by a result of Ball \cite{Ba0} (see Proposition 3.4 in \cite{Pa-sb}). 
\bigbreak

\noindent {\bf Remark.} In the case of projections the connection between comparison and distance inequalities
was established by Ball \cite{Ba2, Ba3}. In particular, it was proved in \cite[Example 2]{Ba3} that
\begin{equation}\label{distance-proj}
c\sqrt{n}\le \max_L d_{\rm vr}(L,\Pi_n)\le \sqrt{n},
\end{equation}
where $\Pi_n$ is the class of projection bodies (origin-symmetric zonoids; see \cite{S2}) in $\R^n,$ $c$ is an absolute constant, and maximum is taken over all origin-symmetric convex bodies in $\R^n.$

The distance $d_{\rm vr}(L,\Pi_n)$ has been studied by several authors. It was introduced in \cite{Ba3} and was proved to be equivalent to the {\it weak-right-hand-Gordon-Lewis} constant of $L$. Also it was connected to the random unconditional constant of the dual space (see Theorem 5 and Proposition 6 in \cite{Ba3}). In \cite{GMP} this distance was called zonoid ratio, and it was proved that it is bounded above by the projection constant of the space. In the same paper  the zonoid ratio was computed for several classical spaces. We refer the interested reader to \cite{GMP}, \cite{GJN}, \cite{GJ}, \cite{GLSW} for more information.

\end{document}